\documentclass[12pt,oneside]{amsart}

\usepackage[english]{babel}
\usepackage{cite}
\usepackage{latexsym}
\usepackage{amssymb,amsthm,amsmath}
\usepackage{hyperref}
\usepackage{graphicx}
\usepackage{longtable}
\usepackage{xcolor}
\usepackage{enumerate}
\usepackage{mathdots}
\usepackage{comment}

\newtheorem{theorem}{Theorem}[section]
\newtheorem*{theorem*}{Theorem}

\newtheorem{remark}[theorem]{Remark}
\newtheorem{proposition}[theorem]{Proposition}

\title{Planar Polynomials arising from Linearized polynomials}

\author{Daniele Bartoli}
\address{Department of Mathematics and Informatics, University of Perugia, Perugia,  Italy}
\email {daniele.bartoli@unipg.it}

\author{Matteo Bonini}
\address{Department of Mathematics, University of Trento, Trento,  Italy}
\email {matteo.bonini@unitn.it}

\date{}

\begin{document}

\maketitle

\begin{abstract}
In this paper we construct planar polynomials of the type $f_{A,B}(x)=x(x^{q^2}+Ax^{q}+Bx)\in \mathbb{F}_{q^3}[x]$, with $A,B \in \mathbb{F}_{q}$. In particular we completely classify the pairs $(A,B)\in \mathbb{F}_{q}^2$ such that $f_{A,B}(x)$ is planar using connections with algebraic curves over finite fields. 	
		
\end{abstract}

	
\section{Introduction}
	Let $q$ be a power of a prime, $n \geq  1$, $\mathbb{F}_{q^n}$ be the finite field with $q^n$ elements and $\mathbb{F}_{q^n}^*= \mathbb{F}_{q^n} \setminus\{0\}$.
A function $f :\mathbb{F}_{q^n} \to \mathbb{F}_{q^n}$ is called planar if for each $a \in \mathbb{F}_{q^n}^*$, the difference
function 
$$D_{f,a} \ :\ \mathbb{F}_{q^n} \to \mathbb{F}_{q^n},\quad  D_{f,a}(x) = f (x+a) - f (x)$$
is a permutation of $\mathbb{F}_{q^n}$. Functions satisfying the previous condition exist only for $q$ odd, since for $q$ even $D_{f,a}(x)=D_{f,a}(x+a)$ holds. In the even characteristic case, functions $f$ such that $D_{f,a}$ is a $2$-to-$1$ map are called  almost perfect nonlinear (APN) and they are connected with  the construction of S-boxes in block ciphers \cite{N1992,N1994}. In the odd characteristic case, planar functions, also called perfect nonlinear functions (or 1-differentiable), are never bijections, since the corresponding derivative 
$D_{f,a}$ is bijective and there exists a unique $x \in \mathbb{F}_{q^n}$ such that $f(x+a)-f(x) = 0$. Planar functions in odd characteristic achieve the best possible differential properties, which makes them useful in the construction of cryptographic protocols, in particular for the design of S-boxes in block chypers; see \cite{N1992,N1994}. Such functions have close connections with projective planes and have been investigated since the seminal paper \cite{DO1968} where the authors showed that, considering two groups $G$ and $H$ of order $n$, every planar function $f$ from $G$ into $H$ gives rise to an affine plane $S(G, H, f)$ of order $n$. Finally, these functions are  strictly related to the construction of relative difference sets \cite{GS1975}, optimal constant-composition codes \cite{DY2005}, secret sharing schemes arising from certain linear codes \cite{CDY2005}, signal sets with good correlation properties \cite{DY2007}, finite semifields \cite{CH2008}. 
More recently, Zhou \cite{25} defined a natural analogue of planar functions for even characteristic: If $q$ is even, a function $f : \mathbb{F}_q \to \mathbb{F}_q$ is planar if, for each nonzero $a \in \mathbb{F}_q$, the function $x\to  f(x + a) + f(x) + ax$ is a permutation of $\mathbb{F}_q$. As shown by Zhou \cite{25} and Schmidt and Zhou \cite{24}, such planar functions have similar properties and applications as their odd characteristic counterparts.

In the past years, many papers have been devoted to existence and non-existence results for planar mappings, using a variety of methods; see
\cite{CM1997,HKNP2008,HS2009,KO2012,CS2016,BN2015,P2016}.

A fundamental tool in our investigations is the Hasse-Weil bound.
\begin{theorem}
Let $\mathcal{X}$ be an absolutely irreducible projective curve defined over $\mathbb{F}_q$. Then the number $\mathcal{X}(\mathbb{F}_q)$ of $\mathbb{F}_q$-rational points of $\mathcal{X}$ satisfies 
\[
||\mathcal{X}(\mathbb{F}_q)|-q-1)|\le  (d-1)(d-2)\sqrt{q},
\]
where $d$ is the degree of $\mathcal{X}$.
\end{theorem}

In this paper we investigate the planarity of $xL(x)$, where $L$ is a linearized polynomial. In this direction, the following results provides a necessary condition.

\begin{proposition}\cite[Proposition 1]{KO2012}
	Let $L_1, L_2 : \mathbb{F}_{p^n} \to \mathbb{F}_{p^n}$ be $\mathbb{F}_q$-linear mappings. If the mapping $L_1 \cdot L_2 : \mathbb{F}_{p^n} \to \mathbb{F}_{p^n}$ is planar then necessarily the mappings $L_1$ and $L_2$ are bijective on $\mathbb{F}_{p^n}$.
\end{proposition}

Polynomial of the type $xL(x)$ belong to the family of so called Dembowski-Ostrom polynomials, that is polynomials of the type 
\begin{equation}\label{DOPolynomial}
\sum_{i,j} a_{ij}X^{p^i+p^j},
\end{equation}
see \cite{DO1968}. It has been conjectured for many years that all planar polynomials belong to such family. Counterexamples have been provided in \cite{CM1997}, where the authors prove that the monomial $X^{(3^{\alpha}+1)/2}$ is planar over $\mathbb{F}_{3^e}$ if and only if $GCD(\alpha,2e)=1$. In the monomial case, the  Dembowski-Ostrom conjecture remains still open for $p\geq 5$. 

There are different possibile definitions of equivalence between planar functions. Two functions $F$ and $F^{\prime}$ from $\mathbb{F}_q$ to itself are called extended affine equivalent (EA-equivalent) if $F= A_1 \circ F \circ A_2+A$, where the mappings $A,A_1,A_2$ are affine, and  $A_1$, $A_2$ are permutations of $\mathbb{F}_q$. The following theorem lists the currently known families of EA-inequivalent planar functions.

\begin{theorem}\label{Th:InequivalentPlanarFunctions}
Let $p$ be an odd prime number. The following are the currently known EA-inequivalent planar functions.
\begin{enumerate}
\item  $x^2$ in $\mathbb{F}_{p^n}$;
\item $x^{p^k+1}$ in $\mathbb{F}_{p^n}$, $k\leq n/2$ and $n/GCD(k,n)$ odd \cite{CM1997,DO1968};
\item $x^{10} + x^6 - x^2$ in $\mathbb{F}_{3^n}$, $n\geq 5$ odd \cite{CM1997};
\item $x^{10} - x^6 - x^2$ in $\mathbb{F}_{3^n}$, $n\geq 5$ odd \cite{DY2006};
\item $x^{p^s+1} -u^{p^k-1}x^{p^k+p^{2k+s}}$ in $\mathbb{F}_{p^{3k}}$, where $GCD(k, 3) = 1$, $k \equiv s \pmod 3$, $s\neq  k$ and $k/GCD(k, s)$ odd,
and $u$ a primitive element of $\mathbb{F}_{p^{3k}}$ \cite[Theorem 1]{ZKW2009};
\item $x^{(3^k+1)/2}$ in $\mathbb{F}_{3^n}$, $k \geq  3$  odd and $GCD(k,n) = 1$ \cite{CM1997,HS1997}.

\end{enumerate}
\end{theorem}

In this paper we investigate the planarity of  $xL_{3}(x)$ on $F_{q^3}$, a particular example of Dembowski-Ostrom polynomials. We make use of the study of curves associated with polynomials $xL_{3}(x)$. The main result of this paper is the following.

\begin{theorem*}
		Let $q$ be a power of an odd prime and consider the function $f_{A,B}(x)\in\mathbb{F}_{q^3}[x]$ of the type 
		\begin{equation}\label{Eq:f_AB}
		f_{A,B}(x)=x(x^{q^2}+Ax^q+Bx),
		\end{equation}
		where $A,B\in \mathbb{F}_q$. Then, $f_{A,B}$ is planar if and only if one of the following condition holds 
		
		\begin{itemize}
		\item $A^3 - 2 A B + 1=0$ and  $A^3\neq \pm1$;  
		\item $A=B^2$, $B^3\neq 1$;
		\item $B=0$, $A^3+1\neq 0$. 
		\end{itemize}
\end{theorem*}

	\section{Construction of planar polynomials of type \texorpdfstring{$x(x^{q^2}+Ax^q+Bx)\in \mathbb{F}_{q^3}[x], A,B \in \mathbb{F}_q$}{Lg}}\label{Section:Dim3}
	
	In this section we provide the complete classification of planar polynomial of the type  $f_{A,B}(x)=x(x^{q^2}+Ax^q+Bx)\in \mathbb{F}_{q^3}[x]$, where $A,B\in \mathbb{F}_q$. Let us denote by $\Box_q$ and $\overline{\mathbb{F}_q}$   the set of squares in $\mathbb{F}_q$ and the algebraic closure of $\mathbb{F}_q$, respectively.

	Note that  $f_{A,B}(x)$ is planar if and only if for each $C\in \mathbb{F}_{q^3}^*$ the linearized polynomial
	$$f_{A,B}(x+C)-f_{A,B}(x)-f_{A,B}(C)=C x^{q^2}+AC x^q+(C^{q^2}+AC^q+2BC)x$$
	is a permutation polynomial. This happens if and only if  the following matrix $M_{A,B,C}$
	\begin{equation}\label{Eq:Matrix}
	\left(\begin{array}{ccc}
	AC^q + 2 B C + C^{q^2}&AC&C\\
		C^q&AC^{q^2} + 2BC^q + C&AC^q\\
		A C^{q^2}&C^{q^2}&A C + 2 B C^{q^2} + C^q\\
	\end{array}\right)
	\end{equation}
	has rank $3$ for each choice of $C\in \mathbb{F}_{q^3}^*$. 
	Let $G_{A,B}(X,Y,T)$ be the determinant of 
	\begin{equation*}
	\left(\begin{array}{ccc}
	AY + 2 B X + T&AX&X\\
		Y&AT + 2BY + X&AY\\
		A T&T&A X + 2 B T + Y\\
	\end{array}\right),
	\end{equation*}
	and define $F_{A,B}(X,Y):=G_{A,B}(X,Y,1)$.
	Clearly,
	\begin{equation}\label{Eq:G}
	\det(M_{A,B,C})=G_{A,B}(C,C^q,C^{q^2}).
	\end{equation}
	It is readily seen that 
	 \begin{eqnarray}\label{Eq_F}
		  F_{A,B}(X,Y)&=&2AB(X^3+Y^3+1)+(2A^2B + 4B^2)(X+Y^2+X^2Y)\\\nonumber
		  &&+(4AB^2 + 2B)(X^2+Y+XY^2)+(2A^3 + 8B^3 + 2)XY.
	\end{eqnarray} 
	
	Our main approach involves specific algebraic curves of degree three, whose reducibility is investigated in the following propositions. 
	
	\begin{proposition}\label{CondizioniGenerali}
		Let $\mathcal{C}_{A,B}$ be a curve defined by $F_{A,B}(X,Y)=0$, $A,B\in \mathbb{F}_{q}^*$, with 
		 \begin{eqnarray*}
		  F_{A,B}(X,Y)&=&2AB(X^3+Y^3+1)+(2A^2B + 4B^2)(X+Y^2+X^2Y)\\
		  &&+(4AB^2 + 2B)(X^2+Y+XY^2)+(2A^3 + 8B^3 + 2)XY.
		  \end{eqnarray*}
		If $\mathcal{C}$ has a line of the form $Y-aX-b$ as a component, then one of the following holds
		\begin{enumerate}
		    \item $A-2B+1=0$ or $(A,B)=(1,-1/2)$; or
		    \item $A^3 - 2 A B + 1=0$ and $A^3+1\neq 0$; or
		    \item $A^2+A+1=0$, $B\in \{A^2,-A^2/2\}$ and $-3\in \square_q$; or
		    \item $A = B^2$; or
		    \item $A^2 + 2 AB - A + 4 B^2 + 2 B + 1=0$ and $-3\in \square_q$.
		\end{enumerate}
	\end{proposition}
	\proof
First we deal with the case $Y-aX-b=Y+X+1$. If such a line is a component, by direct computations	either $A - 2B + 1=0$ or $(A,B)=(1,-1/2)$. On the other hand, it easily seen that if $A - 2B + 1=0$ or $(A,B)=(1,-1/2)$ then  $Y+X+1$ is a factor of $F_{A,B}(X,Y)$. 

By direct computations, $F(X,aX+b)\equiv 0$ if and only if 
\begin{eqnarray*}
B(A + 2 B b + b^2)(A b + 1)&=&0,\\
B(A + a)(A a^2 + 2 B a + 1)&=&0,\\
A^3 b + A^2 B a + A^2B b^2 + 4 A B^2 a b + 2 A B^2 + 3 A B a b^2&&\\
+ 4 B^3 b + 2 B^2 a + 2 B^2 b^2 + 2 B a b + B + b&=&0,\\
A^3a + 2 A^2 B a b + A^2 B + 2AB^2 a^2 + 2 A B^2 b + 3 A B a^2b&&\\
+ 4B^3 a + 4 B^2 a b + 2 B^2 + B a^2 + B b + a&=&0.\\            
\end{eqnarray*}
We distinguish two cases.
\begin{itemize}
        \item $b=-1/A$. 
        \begin{itemize}
        \item If $a=-A$ then $A^3 - 2 A B + 1=0$ or  $A^2B + A - 2B^2=A^2 - 2AB^2 + B=0$. In this last case, either $(A,B)\in \{(1,1),(1,-1/2)\}$ (already done) or $A^2+A+1=0$ (and then $-3\in \square_q$) and $B=A^2,-A^2/2$. 
        \item If $A a^2 + 2 B a + 1=0$ then $A = B^2$ or $A^3 - 2AB + 1=0$. In this last case $A^3+1\neq 0$ otherwise $B=0$. 
        \end{itemize}
        \item  $A + 2 B b + b^2=0$. 
        \begin{itemize}
        \item If $a=-A$ then $A^3 - 2 A B + 1=0$ or $A=B^2$. 
        \item If $A a^2 + 2 B a + 1=0$ then $A - 2B + 1=0$ (already done), or $A=B^2$, or $A^2 + 2 AB - A + 4 B^2 + 2 B + 1=0$ or $A^3 - 2 A B + 1=0$. Note that $A^2 + 2 AB - A + 4 B^2 + 2 B + 1=0$ implies $B=(-A-1+\alpha (A-1))/4$, with $\alpha^2=-3$. Since $A=1$ would imply  $(A,B)\in \{(1,1),(1,-1/2)\}$ (already done) then $A\neq 1$ and  $-3\in \square_q$.
        \end{itemize}
\end{itemize}
	\endproof
	
\begin{proposition}\label{Prop:Factorizations}
Suppose that $\mathcal{C}_{A,B}$, $AB\neq 0$, contains a line of the form $Y-aX-b$ as a component with $(a,b)\neq (-1,-1)$. Then one of the following cases occurs.
\begin{enumerate}
\item $A - 2B + 1=0$ or $(A,B)=(1,-1/2)$  and
$$(Y+X+1)\mid F_{A,B}(X,Y).$$
\item $A^3 - 2 A B + 1=0$, $A^3+1\neq 0$ and 
$$F_{A,B}(X,Y)=(A^2 + AX + Y)(A^2Y + A + X)(A^2X + AY + 1).$$
\item $A=B^2$ and 
$$F_{A,B}(X,Y)=(B^2 + BY + X)(B^2Y + BX + 1)(B^2X + B + Y).$$
\item $-3\in \square_q$, $A^2 + 2 AB - A + 4 B^2 + 2 B + 1=0$ and 
$$(2X- Y-1+\alpha(Y-1)) \mid F_{A,B}(X,Y),$$
where $ \alpha^2=-3$.
\item $-3\in \square_q$, $A^2+A+1=0$, $B\in \{A^2,-A^2/2\}$, and 
$$(2X- Y-1+\alpha(Y-1)) \mid F_{A,B}(X,Y),$$
where $ \alpha^2=-3$.
\end{enumerate}

\end{proposition}	
\proof
Proposition \ref{CondizioniGenerali} applies here. The first three cases are just easy computations. 

Now suppose that $A^2 + 2 AB - A + 4 B^2 + 2 B + 1=0$ and $-3\in \square_q$. Then $B=\frac{-A-1+\alpha(A-1)}{4}$, where $\alpha^2+3=0$ and $F_{A,B}\left(\frac{Y+1+\alpha(Y-1)}{2},Y\right)$ vanishes. 

If  $-3\in \square_q$, $A^2+A+1=0$, and  $B\in \{A^2,-A^2/2\}$ then $A=(-1+\alpha)/2$, where $\alpha^2+3=0$. It is easily seen that $F_{A,B}\left(\frac{Y+1+\alpha(Y-1)}{2},Y\right)$ vanishes.
\endproof

\begin{proposition}\label{Prop:Factorizations_AB0}
The curve $\mathcal{C}_{A,0}$, decomposes in $(A^3+1)XY=0$. The curve $\mathcal{C}_{0,B}$ is reducible only if $B^3+1/8=0$. In this case $(2BX+4B^2Y+1)\mid F_{0,B}(X,Y)$.
\end{proposition}
\proof
If $B=0$, then there is nothing to prove. Suppose now $A=0$ and then $B\neq 0$. Then $F_{0,B}(X,Y)$ reads  
 $$2B(2BY+X)XY + 2BY^2 + (8B^3+2)XY + 4B^2X^2 +  2B(2BY+X).$$
If $F_{0,B}(X,Y)$ is not absolutely irreducible then it contains a line of type $X=a$, $Y=a$, or $X+2BY+a=0$ for some $a\in \overline{\mathbb{F}_q}$. By direct computations, the first two cases imply $B=0$, a  contradiction. The case $X+2BY+a=0$ gives $B^3+1/8=0$ and $a=\frac{1}{2B}$. 
\endproof 

Finally, we need this easy results on linearized polynomials. 
\begin{proposition}\label{Prop:LinearizedPolynomials}
A linearized polynomial $\alpha X^{q^2}+\beta X^{q}+\gamma X$, with $\alpha,\beta,\gamma \in \mathbb{F}_q$ has nonzero roots in $\mathbb{F}_{q^3}$ if and only if $\alpha^3+\beta^3+\gamma^3-3\alpha\beta\gamma= 0$. 
\end{proposition}
\proof
The polynomial $\alpha X^{q^2}+\beta X^{q}+\gamma X$ has a nonzero root if and only if the rank of the matrix 
\begin{equation}\label{Eq:Matrix2}
	\left(\begin{array}{ccc}
	\gamma&\beta&\alpha\\
	\alpha&\gamma&\beta\\
	\beta&\alpha&\gamma\\
	\end{array}\right)
\end{equation}
is smaller than three, that is, $\alpha^3+\beta^3+\gamma^3-3\alpha\beta\gamma= 0$.
\endproof
	
	We are now in a position to prove our main results concerning planar polynomials of type $f_{A,B}(x)$ of $\mathbb{F}_{q^3}$ as in  \eqref{Eq:f_AB}.
	
	\begin{theorem}\label{Th:Main}
		Let $q$ be a power of an odd prime and consider the function $f_{A,B}(x)\in\mathbb{F}_{q^3}[x]$ of the type 
		\[
		f_{A,B}(x)=x(x^{q^2}+Ax^q+Bx),
		\]
		where $A,B\in \mathbb{F}_q$. Then, $f_{A,B}$ is planar if and only if one of the following condition holds 
		
		\begin{itemize}
		\item $A^3 - 2 A B + 1=0$,  $A^3\neq \pm1$;  
		\item $A=B^2$, $B^3\neq 1$;
		\item $B=0$, $A^3+1\neq 0$. 
		\end{itemize}
	\end{theorem}
	
	\proof
	Recall that $f_{A,B}(x)$ is planar if and only if $\det(M_{A,B,C})$ has no roots $C\in \mathbb{F}_{q^3}^*$. 
	We consider first the case $AB\neq 0$. 
	\begin{enumerate}
	\item If $A^3 - 2 A B + 1=0$ and  $A^3+1\neq 0$ then $F_{A,B}(X,Y)=(A^2 + AX + Y)(A^2Y + A + X)(A^2X + AY + 1)$ by Proposition \ref{Prop:Factorizations}.  Therefore 
	$$\det(M_{A,B,C})=(A^2 C^{q^2} + AC + C^q)(A^2C^q + AC^{q^2} + C)(A^2C + AC^q + C^{q^2}).$$
	By Proposition \ref{Prop:LinearizedPolynomials}, all factors have no  roots $C\in \mathbb{F}_{q^3}^*$ if and only if $A^3-1\neq 0$. 
	\item If $A=B^2$ then $F_{A,B}(X,Y)=(B^2 + BY + X)(B^2Y + BX + 1)(B^2X + B + Y)$ by Proposition \ref{Prop:Factorizations}.  Therefore 
	$$\det(M_{A,B,C})=(B^2C^{q^2} + BC^{q} + C)(B^2C^q + BC + C^{q^2})(B^2C + BC^{q^2} + C^q).$$
	By Proposition \ref{Prop:LinearizedPolynomials}, all  factors has no roots $C\in \mathbb{F}_{q^3}^*$ if and only if $B^3-1\neq 0$. 
	\end{enumerate}
	By Propositions \ref{CondizioniGenerali} and \ref{Prop:Factorizations} if  $F_{A,B}(X,Y)$ has a linear factor but $A^3 - 2 A B + 1\neq 0$ and  $A\neq B^2$ then such a factor corresponds to a linearized factor of $\det(M_{A,B,C})$ which is either $C+C^q+C^{q^2}$ or $2C- C^q-C^{q^2}+\alpha(C^q-C^{q^2})$, with $\alpha^2+3=0$ and $-3\in \square_q$. It can be easily verified, using Proposition \ref{Prop:LinearizedPolynomials}, that such a factor has a nonzero root in $\mathbb{F}_{q^3}$ and therefore $f_{A,B}(x)$ is not planar. 
	
	We deal now with the case $AB=0$. 
	\begin{enumerate}
	\item If $B=0$, by Proposition \ref{Prop:Factorizations_AB0} $F_{A,B}(X,Y)=(A^3+1)XY$, that is  $\det(M_{A,B,C})=(A^3+1)C^{q+1}$ and its unique root is $C=0$ and thus $f_{A,B}(x)$ is planar. 
	\item Consider now the case $A=0$ and $B\neq 0$. By Proposition \ref{Prop:Factorizations_AB0}, if $F_{A,B}(X,Y)$ is reducible then $B^3+1/8=0$ and  $(2BX+4B^2Y+1)$ is a factor. Such a factor corresponds to a factor $2BC+4B^2C^q+C^{q^2}$ in $\det(M_{A,B,C})$ and by Proposition \ref{Prop:LinearizedPolynomials} it has nonzero roots $C\in \mathbb{F}_{q^3}$ and so $f_{A,B}(x)$ is not planar. 
	\end{enumerate}
	We are left with the cases in which $F_{A,B}(X,Y)$ is absolutely irreducible. Now consider the polynomial $G_{A,B}(X,Y,T)$ as in \eqref{Eq:G}. Let us fix a normal basis $\{\xi,\xi^q,\xi^{q^2}\}$ of $\mathbb{F}_{q^3}$ over $\mathbb{F}_q$. Then each element $C\in \mathbb{F}_{q^3}^*$ can be written in a unique way as $C=X \xi +Y\xi^q+T\xi^{q^2}$, with $(X,Y,T)\neq (0,0,0)$, $X,Y,T\in \mathbb{F}_q$. The curves defined by $H_{A,B}(X,Y,T)=G_{A,B}(X \xi +Y\xi^q+T\xi^{q^2},X \xi^q +Y\xi^{q^2}+T\xi,X \xi^{q^2} +Y\xi+T\xi^{q})\in \mathbb{F}_{q}[X,Y,Z]$ and $G_{A,B}(X,Y,T)$ are $\mathbb{F}_{q^3}$-isomorphic. Since $F_{A,B}(X,Y)$ and therefore $G_{A,B}(X,Y,T)$ are absolutely irreducible, so it is $H_{A,B}(X,Y,T)$. Since $q>3$, by Hasse-Weil Theorem the curve defined by $H_{A,B}(X,Y,T)=0$ has $\mathbb{F}_{q}$-rational points $(x_0,y_0,t_0)$ which correspond to nonzero roots $C=x_0 \xi +y_0\xi^q+t_0\xi^{q^2}$ of $\det(M_{A,B,C})$. Thus, in all these cases $f_{A,B}(x)$ is not planar.
	\endproof
	
\begin{remark}
By Theorem \ref{Th:Main} the number of pairs $(A,B)\in \mathbb{F}_{q}$ such that $f_{A,B}(x)$ is a planar polynomial in $\mathbb{F}_{q^3}$ is $3q-2-4GCD(3,q-1)$.
\end{remark}

\section{Acknowledgments}
The research was supported by the Italian National Group for Algebraic and Geometric Structures and their Applications (GNSAGA - INdAM).


\end{document}